\documentclass[11pt,a4paper,reqno]{article}
\usepackage{amsmath}
\usepackage{amsthm}
\usepackage{enumerate}
\usepackage{amssymb}
\usepackage{mathrsfs}

\usepackage{verbatim}
\usepackage{url}
\usepackage[latin1]{inputenc}

\usepackage{caption}
\usepackage{graphicx,color}
\def\Z{{\mathbb Z}}

\def\R{{\mathbb R}}

\def\coex{{\rm Coex}}
\def\gex{{\rm Geos}}
\def\tree{{\mathscr T}}
\def\ball{{\rm Ball}}
\def\geo{{\rm geo}}

\def\be{\begin{equation}}
\def\ee{\end{equation}}
\def\bea{\begin{equation*}}
\def\eea{\end{equation*}}
\def\bal{\begin{aligned}}
\def\eal{\end{aligned}}

\def\eps{\varepsilon}

\def\Pr{{\mathbb P}}
\DeclareMathOperator{\E}{{\mathbb E}}

\newtheorem{thm}{Theorem}

\newtheorem{cor}[thm]{Corollary}

\theoremstyle{remark}

\newtheorem{preex}[thm]{Example}

\theoremstyle{definition}

\begin{document}

\title{Existence and coexistence in first-passage percolation}
\date{}
\author{Daniel Ahlberg\thanks{Department of Mathematics, Stockholm University.}}
\maketitle

\begin{abstract}
We consider first-passage percolation with i.i.d.\ non-negative weights coming from some continuous distribution under a moment condition. We review recent results in the study of geodesics in first-passage percolation and study their implications for the multi-type Richardson model. In two dimensions this establishes a dual relation between the existence of infinite geodesics and coexistence among competing types. The argument amounts to making precise the heuristic that infinite geodesics can be thought of as `highways to infinity'. We explain the limitations of the current techniques by presenting a partial result in dimensions $d>2$.
\end{abstract}

\section{Introduction}

In first-passage percolation the edges of the $\Z^d$ nearest neighbour lattice, for some $d\ge2$, are equipped with non-negative i.i.d.\ random weights $\omega_e$, inducing a random metric $T$ on $\Z^2$ as follows: For $x,y\in\Z^d$, let
\begin{equation}\label{eq:Tdef}
\textstyle{T(x,y):=\inf\big\{\sum_{e\in\pi}\omega_e:\pi\text{ is a self-avoiding path from $x$ to $y$}\big\}.}
\end{equation}
Since its introduction in the 1960s, by Hammersley and Welsh~\cite{hamwel65}, a vast body of literature has been generated seeking to understand the large scale behaviour of distances, balls and geodesics in this random metric space. The state of the art has been summarized in various volumes over the years, including~\cite{aufdamhan17,howard04,kesten86,smywie78}. We will here address questions related to geodesics, and shall for this reason make the common assumption that the edge weights are sampled from a continuous distribution. Since many of the results we shall rely on require a moment condition for their conclusions to hold, we shall assume in what follows that $\E[Y^d]<\infty$, where $Y$ denotes the minimum weight among the $2d$ edges connected to the origin.

In the 1960s, the study of first-passage percolation led to the development of an ergodic theory for subadditive ergodic sequences, culminating with the ergodic theorem due to Kingman~\cite{kingman68}. As a consequence thereof, one obtains the existence of a norm $\mu:\R^d\to[0,\infty)$, simply referred to as the \emph{time constant}, such that for every $z\in\Z^d$, almost surely,
$$
\lim_{n\to\infty}\frac1nT(0,nz)=\mu(z).
$$
Richardson~\cite{richardson73}, and later work of Cox and Durrett~\cite{coxdur81}, extended the above \emph{radial} convergence to \emph{simultaneous} convergence in all directions. Their results show that the ball $\{z\in\Z^d:T(0,z)\le t\}$ in the metric $T$ once rescaled by $1/t$ approaches the unit ball in the norm $\mu$. The unit ball in $\mu$, henceforth denoted by $\ball:=\{x\in\R^d:\mu(x)\le1\}$, is therefore commonly referred to as the \emph{asymptotic shape}, and known to be compact and convex with non-empty interior. In addition, the shape retains the symmetries of $\Z^d$. However, little else is known regarding the properties of the shape in general. This, we shall see, is a major obstacle for our understanding of several other features of the model.

Although questions regarding geodesics were considered in the early work of Hammersley and Welsh, it took until the mid 1990s before Newman~\cite{newman95} together with his co-authors~\cite{licnew96,licnewpiz96,newpiz95} initiated a systematic study of the geometry of geodesics in first-passage percolation. Under the assumption of continuous weights there is almost surely a unique path attaining the minimum in~\eqref{eq:Tdef}; we shall denote this path $\geo(x,y)$ and refer to it as the \emph{geodesic} between $x$ and $y$. The graph consisting of all edges on $\geo(0,y)$ for some $y\in\Z^d$ is a tree spanning the lattice. Understanding the properties of this object, such as the number of topological ends, leads one to the study of \emph{infinite} geodesics, i.e.\ infinite paths $g=(v_1,v_2,\ldots)$ of which every finite segment is a geodesic.
We shall write $\tree_0$ for the collection of infinite geodesics starting at the origin. A simple compactness argument shows that the cardinality $|\tree_0|$ of $\tree_0$ is always at least one.\footnote{Consider the sequence of finite geodesics between the origin and $n{\bf e}_1$, where ${\bf e}_1$ denotes the first coordinate vector. Since the number of edges that connect to the origin is finite, one of them must be traversed for infinitely many $n$. Repeating the argument results in an infinite path which by construction is a geodesic.}
In two dimensions, Newman~\cite{newman95} predicted that $|\tree_0|=\infty$ almost surely, and proved this under an additional assumption of uniform curvature of the asymptotic shape, which remains unverified to this day.

As a means to make rigorous progress on Newman's prediction, H\"aggstr\"om and Pemantle~\cite{hagpem98} introduced a model for competing growth on $\Z^d$, for $d\ge2$, known as the \emph{two-type Richardson model}. In this model, two sites $x$ and $y$ are initially coloured red and blue respectively. As time evolves an uncoloured site turns red at rate $1$ times the number of red neighbours, and blue at rate $\lambda$ times the number or blue neighbours. A central question of interest is for which values of $\lambda$ there is positive probability for both colours to coexist, in the sense that they both are responsible for the colouring of infinitely many sites.

There is an intimate relation between the existence of infinite geodesics and coexistence in the Richardson model that we shall pay special interest in. In the case of equal strength competitors ($\lambda=1$), one way to construct the two-type Richardson model is to equip the edges of the $\Z^d$ lattice with independent exponential weights, thus exhibiting a direct connection to first-passage percolation.
The set of sites eventually coloured red in the two-type Richardson model is then equivalent to the set of sites closer to $x$ than $y$ in the first-passage metric. That is, an analogous way to phrase the question of coexistence is whether there are infinitely many points closer to $x$ than $y$ as well as infinitely many points closer to $y$ than $x$ in the first-passage metric. As before, a compactness argument will show that on the event of coexistence there are disjoint infinite geodesics $g$ and $g'$ that respectively originate from $x$ and $y$. H\"aggstr\"om and Pemantle~\cite{hagpem98} showed that, for $d=2$, coexistence of the two types occurs with positive probability, and deduced as a corollary that
$$
\Pr(|\tree_0|\ge2)>0.
$$
Their results were later extended to higher dimensions and more general edge weight distributions in parallel by Garet and Marchand~\cite{garmar05} and Hoffman~\cite{hoffman05}. In a later paper, Hoffman~\cite{hoffman08} showed that in two dimensions coexistence of four different types has positive probability, and that $\Pr(|\tree_0|\ge4)>0$. The best currently known general lower bound on the number of geodesics is a strengthening of Hoffman's result due to Damron and Hanson~\cite{damhan14}, showing that
$$
\Pr(|\tree_0|\ge4)=1.
$$

In this paper we shall take a closer look at the relation between existence of infinite geodesics and coexistence in competing first-passage percolation. We saw above that on the event of coexistence of various types, a compactness argument gives the existence of equally many infinite geodesics. It is furthermore conceivable that it is possible to locally modify the edge weight in such a way that these geodesics are re-routed through the origin. Conversely, interpreting infinite geodesics as `highways to infinity', along which the different types should be able to escape their competitors, it seems that the existence of a given number of geodesics should accommodate an equal number of surviving types. These heuristic arguments suggest a duality between existence and coexistence, and it is this dual relation we shall make precise.

Given sites $x_1,x_2,\ldots,x_k$ in $\Z^d$, we let $\coex(x_1,x_2,\ldots,x_k)$ denote the event that for every $i=1,2,\ldots,k$ there are infinitely many sites $z\in\Z^d$ for which the distance $T(x_j,z)$ is minimized by $j=i$. (The continuous weight distribution assures that there are almost surely no ties.) In two dimensions the duality between existence and coexistence that we prove takes the form:
\begin{equation}\label{eq:duality}
\exists\, x_1,x_2,\ldots,x_k\text{ such that }\Pr\big(\coex(x_1,x_2,\ldots,x_k)\big)>0\quad\Leftrightarrow\quad\Pr(|\tree_0|\ge k)>0.
\end{equation}
Turning the above heuristic into a proof is more demanding that it may seem. In order to derive the relation in~\eqref{eq:duality} we shall rely on the recently developed ergodic theory for infinite geodesics. This theory has its origins in the work of Hoffman~\cite{hoffman05,hoffman08}, and was developed further by Damron and Hanson~\cite{damhan14,damhan17}, before it reached its current status in work of Ahlberg and Hoffman~\cite{ahlhof}. The full force of this theory is currently restricted to two dimensions, which prevents us from obtaining an analogue to~\eqref{eq:duality} in higher dimensions. In higher dimensions we deduce a partial result based on results of Damron and Hanson~\cite{damhan14} and Nakajima~\cite{nakajima}.

\subsection{The dual relation}

Before we state our results formally, we remind the reader that $Y$ denotes the minimum weight among the $2d$ edges connected to the origin. We recall (from~\cite{coxdur81}) that $\E[Y^d]<\infty$ is both necessary and sufficient in order for the shape theorem to hold in dimension $d\ge2$.

\begin{thm}\label{thm:duality}
Consider first-passage percolation on the square lattice with continuous edge weights satisfying $\E[Y^2]<\infty$. For any $k\ge1$, including $k=\infty$, and $\eps>0$ we have:
\begin{enumerate}[\quad (i)]
\item If $\Pr\big(\coex(x_1,x_2,\ldots,x_k)\big)>0$ for some $x_1,x_2,\ldots,x_k$ in $\Z^2$, then $\Pr(|\tree_0|\ge k)=1$.
\item If $\Pr(|\tree_0|\ge k)>0$, then $\Pr\big(\coex(x_1,x_2,\ldots,x_k)\big)>1-\eps$ for some $x_1,x_2,\ldots,x_k$ in $\Z^2$.
\end{enumerate}
\end{thm}

In dimensions higher than two we shall establish parts of the above dual relation, and recall next some basic geometric concepts in order to state this result precisely. A hyperplane in the $d$-dimensional Euclidean space divides $\R^d$ into two open half-spaces. A \emph{supporting hyperplane} to a convex set $S\subset\R^d$ is a hyperplane that contains some boundary point of $S$ and contains all interior points of $S$ in one of the two half-spaces associated to the hyperplane. It is well-known that for every boundary point of a convex set $S$ there exists a supporting hyperplane that contains that point. A supporting hyperplane to $S$ is called a \emph{tangent hyperplane} if it is the unique supporting hyperplane containing some boundary point of $S$. Finally, we define the number of \emph{sides} of a compact convex set $S$ as the number of (distinct) tangent hyperplanes to $S$. Hence, the number of sides is finite if and only if $S$ is a (finite) convex polygon ($d=2$) or convex polytope ($d\ge3$). A deeper account on convex analysis can be found in~\cite{rockafellar70}.

\begin{thm}\label{thm:highdim}
Consider first-passage percolation on the $d$-dimensional cubic lattice, for $d\ge2$, with continuous edge weights. For any $k\ge1$, including $k=\infty$, and $\eps>0$ we have
\begin{enumerate}[\quad (i)]
\item If $\E[\exp(\alpha\omega_e)]<\infty$ and $\Pr\big(\coex(x_1,x_2,\ldots,x_k)\big)>0$ for some $\alpha>0$ and $x_1,x_2,\ldots,x_k$ in $\Z^d$, then $\Pr(|\tree_0|\ge k)=1$.
\item If $\E[Y^d]<\infty$ and $\ball$ has at least $k$ sides, then $\Pr\big(\coex(x_1,x_2,\ldots,x_k)\big)>1-\eps$ for some $x_1,x_2,\ldots,x_k$ in $\Z^d$.
\end{enumerate}
\end{thm}

In Section~\ref{sec:geodesics} we shall review the recent development in the study of infinite geodesics that will be essential for the deduction, in Section~\ref{sec:duality}, of the announced dual result. Finally, in Section~\ref{sec:highdim}, we prove the partial result in higher dimensions.

\subsection{A mention of our methods}

One aspect of the connection between existence and coexistence is an easy observation, and was hinted at already above. Namely, if $\gex(x_1,x_2,\ldots,x_k)$ denotes the event that there exist $k$ pairwise disjoint infinite geodesics, each originating from one of the points $x_1,x_2,\ldots,x_k$, then
\begin{equation}\label{eq:inclusion}
\coex(x_1,x_2,\ldots,x_k)\subseteq\gex(x_1,x_2,\ldots,x_k).
\end{equation}
To see this, let $V_i$ denote the set of sites closer to $x_i$ than to any other $x_j$, for $j\neq i$, in the first-passage metric. (Note that $T(x,y)\neq T(z,y)$ for all $x,y,z\in\Z^2$ almost surely, due to the assumptions of continuous weights.\footnote{This will be referred to as having \emph{unique passage times}.}) On the event $\coex(x_1,x_2,\ldots,x_k)$ each set $V_i$ is infinite, and for each $i$ a compactness argument gives the existence of an infinite path contained in $V_i$, which by construction is a geodesic. Since $V_1,V_2,\ldots,V_k$ are pairwise disjoint, due to uniqueness of geodesics, so are the resulting infinite geodesics.

Let $\mathscr{N}$ denote the maximal number of pairwise disjoint infinite geodesics. Since $\mathscr{N}$ is invariant with respect to translations (and measurable) it follows from the ergodic theorem that $\mathscr{N}$ is almost surely constant. Hence, positive probability for coexistence of $k$ types implies the almost sure existence of $k$ pairwise disjoint geodesics. That $|\tree_0|\le\mathscr{N}$ is trivial, given the tree structure of $\tree_0$. The inequality is in fact an equality, which was established by different means in~\cite{ahlhof,nakajima}. Together with~\eqref{eq:inclusion}, this resolves the first part of Theorems~\ref{thm:duality} and~\ref{thm:highdim}.

Above it was suggested that infinite geodesics should, at least heuristically, be thought of as `highways to infinity' along which the different types may escape the competition. The concept of Busemann functions, and their properties, will be central in order to make this heuristic precise. These functions have their origin in the work of Herbert Busemann~\cite{busemann55} on metric spaces. In first-passage percolation, Busemann-related limits first appeared in the work of Newman~\cite{newman95} as a means to describe the microscopic structure of the boundary (or surface) of a growing ball $\{z\in\Z^2:T(0,z)\le t\}$ in the first-passage metric. Later work of Hoffman~\cite{hoffman05,hoffman08} developed a method to describe asymptotic properties of geodesics via the study of Busemann functions. Hoffman's approach has since become indispensable in the study of various models for spatial growth, including first-passage percolation~\cite{damhan14,damhan17,ahlhof}, the corner growth model~\cite{georassep17a,georassep17b} and random polymers~\cite{georassep16,albrassim20}. In a tangential direction, Bakhtin, Cator and Khanin~\cite{bakcatkha14} used Busemann functions to construct stationary space-time solutions to the one-dimensional Burgers equation, inspired by earlier work of Cator and Pimentel~\cite{catpim12}.

Finally, we remark that (for $d=2$) it is widely believed that the asymptotic shape is not a polygon, in which case it follows from~\cite{hoffman08} that both $\Pr(|\tree_0|=\infty)=1$ and that for every $k\ge1$ there are $x_1,x_2,\ldots,x_k$ such that $\Pr(\coex(x_1,x_2,\ldots,x_k))>0$. The latter was extended to infinite coexistence by Damron and Hochman~\cite{damhoc13}. Thus, proving that the asymptotic shape is non-polygonal would make our main theorem obsolete. However, understanding the asymptotic shape is a notoriously hard problem, which is the reason an approach sidestepping Newman's curvature assumption has been developed in the first place.

\section{Geodesics and Busemann functions}\label{sec:geodesics}

In this section we review the recent developments in the study of infinite geodesics in first-passage percolation. We shall focus on the two-dimensional setting, and remark on higher dimensions only at the end. We make no claim in providing a complete account of previous work, and instead prefer to focus on the results that will be of significance for the purposes of this paper. A more complete description of these results, save those reported in the more recent studies~\cite{ahlhof,nakajima}, can be found in~\cite{aufdamhan17}.

\subsection{Geodesics in Newman's contribution to the 1994 ICM proceedings}

The study of geodesics in first-passage percoalation was pioneered by Newman and co-authors~\cite{licnew96,licnewpiz96,newman95,newpiz95} in the mid 1990s. Their work gave rise to a precise set of predictions for the structure of infinite geodesics. In order to describe these predictions we shall need some notation. First, we say that an infinite geodesic $g=(v_1,v_2,\ldots)$ has \emph{asymptotic direction} $\theta$, in the unit circle $S^1:=\{x\in\R^2:|x|=1\}$, if the limit $\lim_{k\to\infty}v_k/|v_k|$ exists and equals $\theta$. Second, two infinite geodesics $g$ and $g'$ are said to \emph{coalesce} if their symmetrical difference $g\Delta g'$ is finite. The predictions originating from the work of Newman and his collaborators can be summarized as, under mild conditions on the weight distribution, the following should hold:
\begin{enumerate}[\quad\em (a)]
\item with probability one, every infinite geodesic has an asymptotic direction;
\item for every direction $\theta$, there is an almost surely unique geodesic in $\tree_0$ with direction $\theta$;
\item for every direction $\theta$, any two geodesics with direction $\theta$ coalesce almost surely.
\end{enumerate}
In particular, these statements would imply that $|\tree_0|=\infty$ almost surely.

Licea and Newman~\cite{newman95,licnew96} proved conditional versions of these statements under an additional curvature assumption of the asymptotic shape. While this assumption seems plausible for a large family of edge weight distributions, there is no known example for which it has been verified. Rigorous proofs of the corresponding statements for a rotation invariant first-passage-like model, where the asymptotic shape is known to be a Euclidean disc, has been obtained by Howard and Newman~\cite{hownew01}. Since proving properties like strict convexity and differentiability of the boundary of the asymptotic shape in standard first-passage percolation appears to be a major challenge, later work has focused on obtaining results without assumptions on the shape.

\subsection{Busemann functions}

Limits reminiscent of Busemann functions first appeared in the first-passage literature in the work of Newman~\cite{newman95}, as a means of describing the microscopic structure of the boundary of a growing ball in the first passage metric. The method for describing properties of geodesics via Busemann functions developed in later work of Hoffman~\cite{hoffman05,hoffman08}.

Given an infinite geodesic $g=(v_1,v_2,\ldots)$ in $\tree_0$ the \emph{Busemann function} $B_g:\Z^2\times\Z^2\to\R$ of $g$ is defined as the limit
\begin{equation}\label{eq:Busemann}
B_g(x,y):=\lim_{k\to\infty}\big[T(x,v_k)-T(y,v_k)\big].
\end{equation}
As observed by Hoffman~\cite{hoffman05}, with probability one the limit in~\eqref{eq:Busemann} exists for every $g\in\tree_0$ and all $x,y\in\Z^2$, and satisfies the following properties:
\begin{itemize}
\item $B_g(x,y)=B_g(x,z)+B_g(z,y)$ for all $x,y,z\in\Z^2$;
\item $|B_g(x,y)|\le T(x,y)$;
\item $B_g(x,y)=T(x,y)$ for all $x,y\in g$ such that $x\in\geo(0,y)$.
\end{itemize}

In~\cite{hoffman05} Hoffman used Busemann functions to establish that there are at least two disjoint infinte geodesics almost surely. In~\cite{hoffman08} he used Busemann functions to associate certain infinite geodesics with sides (tangent lines) of the asymptotic shape. The approach involving Busemann functions in order to study infinite geodesics was later developed further in work by Damron and Hanson~\cite{damhan14,damhan17} and Ahlberg and Hoffman~\cite{ahlhof}. Studying Busemann functions of geodesics, as opposed to the geodesics themselves, has allowed these authors to establish rigorous versions of Newman's predictions regarding the structure of geodesics. Describing parts of these results in detail will be essential in order to understand the duality between existence of geodesics and coexistence in competing first-passage percolation.


\subsection{Linearity of Busemann functions}

We shall call a linear functional $\rho:\R^2\to\R$ \emph{supporting} if the line $\{x\in\R^2:\rho(x)=1\}$ is a supporting line to $\partial\ball$ through some point, and \emph{tangent} if $\{x\in\R^2:\rho(x)=1\}$ is the unique supporting line (i.e.\ the tangent line) through some point of $\partial\ball$. Given a supporting functional $\rho$ and a geodesic $g\in\tree_0$ we say that the Busemann function of $g$ is \emph{asymptotically linear} to $\rho$ if
\begin{equation}\label{eq:linearity}
\limsup_{|y|\to\infty}\frac{1}{|y|}\big|B_g(0,y)-\rho(y)\big|=0.
\end{equation}

Asymptotic linearity of Busemann functions is closely related to asymptotic directions of geodesics in the sense that~\eqref{eq:linearity}, together with the third of the properties of Busemann functions exhibited by Hoffman, provides information on the direction of $g=(v_1,v_2,\ldots)$: The set of limit points of the sequence $(v_k/|v_k|)_{k\ge1}$ is contained in the arc $\{x\in S^1:\mu(x)=\rho(x)\}$, corresponding to a point or a flat edge of $\partial\ball$.

Building on the work of Hoffman~\cite{hoffman08}, Damron and Hanson~\cite{damhan14} showed that for every tangent line of the asymptotic shape there exists a geodesic whose Busemann function is described by the corresponding linear functional. In a simplified form their result reads as follows:

\begin{thm}\label{thm:DH}
For every tangent functional $\rho:\R^2\to\R^2$ there exists, almost surely, a geodesic in $\tree_0$ whose Busemann function is asymptotically linear to $\rho$.
\end{thm}

While the work of Damron and Hanson proves \emph{existence} of geodesics with linear Busemann functions, later work of Ahlberg and Hoffman~\cite{ahlhof} has established that \emph{every} geodesic has a linear Busemann function, and that the associated linear functionals are \emph{unique}. We summarize these results in the next couple of theorems.

\begin{thm}\label{thm:linearity}
With probability one, for every geodesic $g\in\tree_0$ there exists a supporting functional $\rho:\R^2\to\R$ such that the Busemann function of $g$ is asymptotically linear to $\rho$.
\end{thm}

To address uniqueness, note that the set of supporting functionals is naturally parametrized by the direction of their gradients. Due to convexity of the shape, these functionals stand in 1-1 correspondence with the unit circle $S^1$. We shall from now on identify the set of supporting functionals with $S^1$.

\begin{thm}\label{thm:ergodicity}
There exists a closed (deterministic) set $\mathscr{C}\subseteq S^1$ such that, with probability one, the (random) set of supporting functionals $\rho$ for which there exists a geodesic in $\tree_0$ with Busemann function asymptotically linear to $\rho$ equals $\mathscr{C}$. Moreover, for every $\rho\in\mathscr{C}$ we have
$$
\Pr\big(\exists\text{ two geodesics in $\tree_0$ with Busemann function linear to }\rho\big)=0.
$$
\end{thm}

From Theorem~\ref{thm:DH} it follows that $\mathscr{C}$ contains all tangent functionals. As a consequence, if $\ball$ has at least $k$ sides (i.e.\ tangent lines), then we have $|\tree_0|\ge k$ almost surely. On the other hand, it follows from Theorem~\ref{thm:linearity} that every geodesic has a linear Busemann function, and by Theorem~\ref{thm:ergodicity} that the set of linear functionals describing these Busemann functions is deterministic. Consequently, if with positive probability $\tree_0$ has size at least $k$, then by the uniqueness part of Theorem~\ref{thm:ergodicity} the set $\mathscr{C}$ has cardinality at least $k$, so that there exist $k$ geodesics described by distinct linear functionals almost surely. All these observations will be essential in proving part~\emph{(ii)} of Theorem~\ref{thm:duality}.

Due to the connection between asymptotic directions and linearity of Busemann functions mentioned above, Theorems~\ref{thm:DH}-\ref{thm:ergodicity} may be seen as rigorous, although somewhat weaker, versions of Newman's predictions~\emph{(a)-(b)}. The rigorous results are weaker in the sense that we do not know whether $\mathscr{C}$ equals $S^1$ or not. Note, however, that Theorem~\ref{thm:ergodicity} provides an `ergodic theorem' in this direction. As we shall describe next, the cited papers provide a rigorous version also of~\emph{(c)}.

\subsection{Coalescence}

An aspect of the above development that we have ignored so far is that of coalescence. For instance, Theorem~\ref{thm:DH} is a simplified version of a stronger statement proved in~\cite{damhan14}, namely that for every tangent functional $\rho:\R^2\to\R$ there exists, almost surely, a family of geodesics $\Gamma=\{\gamma_z:z\in\Z^2\}$, where $\gamma_z\in\tree_z$, such that any one geodesic in $\Gamma$ has Busemann function linear to $\rho$ and any two geodesics in $\Gamma$ coalesce. (The latter of course implies that the Busemann functions of all geodesics in $\Gamma$ coincide.) In a similar spirit, we have the following from~\cite{ahlhof}:

\begin{thm}\label{thm:coalescence}
For every supporting functional $\rho\in\mathscr{C}$, with probability one, any two geodesics $g\in\tree_y$ and $g'\in\tree_z$ with Busemann function asymptotically linear to $\rho$ coalesce.
\end{thm}

We remark that coalescence was irrelevant for the proof of Theorem~\ref{thm:DH} in~\cite{damhan14}, but instrumental for the deduction of Theorems~\ref{thm:linearity} and~\ref{thm:ergodicity} in~\cite{ahlhof}. In short, the importance of coalescence lies in the possibility to apply the ergodic theorem to asymptotic properties of shift invariant families of coalescing geodesics, resulting in the ergodic properties of Theorem~\ref{thm:ergodicity}.

The results described above together address the cardinality of the set $\tree_0$. Recall that $\mathscr{N}$ denotes the maximal number of pairwise disjoint infinite geodesics and is almost surely constant. The following was first established in~\cite{ahlhof}, and can be derived as a corollary to Theorems~\ref{thm:linearity}-\ref{thm:coalescence}. A more direct argument, assuming a stronger moment condition, was later given by Nakajima~\cite{nakajima}.

\begin{cor}\label{cor:cardinality}
With probability one $|\tree_0|$ is constant and equal to $\mathscr{N}$.
\end{cor}

To see how the corollary follows, first note that clearly $|\tree_0|\le\mathscr{N}$. In addition, $|\mathscr{C}|\le|\tree_0|$ almost surely due to the ergodic part of Theorem~\ref{thm:ergodicity}, and in the case that $\mathscr{C}$ is finite, equality follows from Theorem~\ref{thm:linearity} and the uniqueness part of Theorem~\ref{thm:ergodicity}. Consequently, also $|\tree_0|$ is almost surely constant. Finally, it follows from the coalescence property in Theorem~\ref{thm:coalescence} that either $|\tree_0|$ (and therefore also $\mathscr{N}$) is almost surely infinite, or $|\tree_0|=\mathscr{N}=k$ holds almost surely for some finite $k$, leading to the claimed result.

\subsection{Geodesics in higher dimensions}

Whether the description of geodesics detailed above remains correct also in higher dimensions is at this point unknown. Although it has been suggested that coalescence should fail for large $d$, it seems plausible that results analogous to Theorems~\ref{thm:DH}-\ref{thm:ergodicity} should hold for all $d\ge2$, and that an analogue to Theorem~\ref{thm:coalescence} could hold for small $d$. See recent work of Alexander~\cite{alexander} for a further discussion of these claims. Indeed, establishing the existence of coalescing families of geodesics in the spirit of~\cite{damhan14} also in three dimensions should be considered a major open problem.

What is known is that the argument behind Theorem~\ref{thm:DH} can be extended to all dimensions $d\ge2$ under minor adjustments; see~\cite{bridamhan}. However, the proofs of Theorems~\ref{thm:linearity}-\ref{thm:coalescence} exploit planarity in a much more fundamental way, and are not known to extend to higher dimensions. On the other hand, an argument of Nakajima~\cite{nakajima} shows that Corollary~\ref{cor:cardinality} remains valid in all dimensions under the additional condition that $\E[\exp(\alpha\omega_e)]<\infty$ for some $\alpha>0$. These properties will be sufficient in order to prove Theorem~\ref{thm:highdim}.

\section{The dual relation in two dimensions}\label{sec:duality}

With the background outlined in the previous section we are now ready to prove Theorem~\ref{thm:duality}. We recall that, with probability one, by Theorem~\ref{thm:linearity} every geodesic has an asymptotically linear Busemann function, and by Theorem~\ref{thm:ergodicity} there is a deterministic set $\mathscr{C}$ of linear functionals that correspond to these Busemann functions. Moreover, for each $\rho\in\mathscr{C}$, by Theorem~\ref{thm:ergodicity} there is for every $z\in\Z^2$ an almost surely unique geodesic in $\tree_z$ with Busemann function asymptotically linear to $\rho$, and by Theorem~\ref{thm:coalescence} these geodesics coalesce almost surely. In particular $|\tree_0|=|\mathscr{C}|$ almost surely, and we shall in the sequel write $B_\rho$ for the Busemann function of the almost surely unique geodesic (in $\tree_0$) corresponding to $\rho$.


\subsection{Part~\emph{(i)}: Coexistence implies existence}

The short proof of part~\emph{(i)} is an easy consequence of Corollary~\ref{cor:cardinality}. Suppose that for some choice of $x_1,x_2,\ldots,x_k$ in $\Z^2$ we have $\Pr(\coex(x_1,x_2,\ldots,x_k))>0$. By~\eqref{eq:inclusion} we have $\Pr(\mathscr{N}\ge k)>0$, and since $\mathscr{N}$ is almost surely constant it follows from Corollary~\ref{cor:cardinality} that
$$
\Pr(|\tree_0|\ge k)=1.
$$

While the above argument is short, it hides much of the intuition for why the implication holds. We shall therefore give a second argument based on coalescence that may be more instructive, even if no more elementary. This argument will make explicit the heuristic that geodesics are `highways to infinity' along which the different types will have to move in order to escape the competition.

Before attending to the proof,
we claim that for any $\rho\in\mathscr{C}$ we have
\begin{equation}\label{eq:nonzero}
\Pr\big(B_\rho(x,y)\neq0\text{ for all }x\neq y\big)=1.
\end{equation}
To see this, let $A_\rho$ denote the event that for each $z$ in $\Z^2$ there is a unique geodesic $g_z$ in $\tree_z$ corresponding to $\rho$, and that all these geodesics coalesce, so that $A_\rho$ has measure one. We note that on the event $A_\rho$ coalescence of the geodesics $\{g_z:z\in\Z^2\}$ implies that for any $x,y\in\Z^2$ the limit $B_\rho(x,y)$ (which is defined through~\eqref{eq:Busemann} for $g=g_0$) is attained after a finite number of steps. More precisely, on the event $A_\rho$, for any $x,y\in\Z^2$ and $v$ contained in $g_x\cap g_y$ we have
$$
B_\rho(x,y)=T(x,v)-T(y,v).
$$
Hence,~\eqref{eq:nonzero} follows due to unique passage times.

We now proceed with the second proof. Again by Corollary~\ref{cor:cardinality}, either $\tree_0$ is almost surely infinite, in which case there is nothing to prove, or $\Pr(|\tree_0|=k)=1$ for some integer $k\ge1$. We shall suppose the latter, and argue that for any choice of $x_1,x_2,\ldots,x_{k+1}$ in $\Z^2$ we have $\Pr(\coex(x_1,x_2,\ldots,x_{k+1}))=0$.

On the event that $\tree_0$ is almost surely finite, $\mathscr{C}$ is in one-to-one correspondence with the elements of $\tree_0$. It  follows from~\eqref{eq:nonzero} that for any $g\in\tree_0$ the Busemann function $B_g(0,x)$ has a unique minimizer over finite subsets of $\Z^2$ almost surely. The last statement can be rephrased in terms of competition between a finite number of types as follows: For each geodesic $g$ in $\tree_0$ there will be precisely one type that reaches infinitely many sites along $g$ almost surely; it is the one whose starting position minimizes $B_g(0,x_i)$. Hence, if $|\tree_0|=k$ almost surely, but there are $k+1$ competing types, then at least one of them will not reach infinitely many sites along any geodesic in $\tree_0$. Suppose that the type left out starts at a site $x$. Since for each geodesic in $\tree_x$ there is a geodesic in $\tree_0$ with which it coalesces (as of Theorem~\ref{thm:coalescence}), it follows that for each geodesic $g\in\tree_x$ the type starting at $x$ will be closer than the other types to at most finitely many sites along $g$. Choose $n$ so that these sites are all within distance $n$ from $x$. Finally, note that for at most finitely many sites $z$ in $\Z^2$ the (finite) geodesic from $x$ to $z$ will diverge from all geodesics in $\tree_x$ within distance $n$ from $x$. Consequently, all but finitely many sites in $\Z^2$ will lie closer to the starting point of some other type, implying that the $k+1$ types do not coexist.

\subsection{Part~\emph{(ii)}: Existence implies coexistence}

Central in the proof of part~\emph{(ii)} is the linearity of Busemann functions.
The argument that follows is a modern take on an argument originally due to Hoffman~\cite{hoffman08}.

Let $k$ be an integer and suppose that $|\tree_0|\ge k$ with positive probability. Then, indeed, $|\tree_0|=|\mathscr{C}|\ge k$ almost surely. Fix $\eps>0$ and let $\rho_1,\rho_2,\ldots,\rho_k$ be distinct elements of $\mathscr{C}$. In order to show that $\Pr(\coex(x_1,x_2,\ldots,x_k))>1-\eps$ for some choice of $x_1,x_2,\ldots,x_k$, we shall choose these points so that with probability $1-\eps$ we have $B_{\rho_i}(x_i,x_j)<0$ for all $i=1,2,\ldots,k$ and $j\neq i$. On this event, for each $i$, the site $x_i$ is closer to all points along the geodesic in $\tree_{x_i}$ corresponding to $\rho_i$ than any of the $x_j$ for $j\neq i$, implying that $\coex(x_1,x_2,\ldots,x_k)$ occurs.

Given $\rho\in\mathscr{C}$, $z\in\Z^2$, $\delta>0$ and $M\ge1$ we let $A_\rho(z,\delta,M)$ denote the event that
$$
\big|B_\rho(z,y)-\rho(y-z)\big|<\delta|y-z|\quad\text{for all }|y-z|\ge M.
$$
Due to linearity of Busemann functions (Theorems~\ref{thm:linearity} and~\ref{thm:ergodicity}) there exists for every $\rho\in\mathscr{C}$ and $\delta,\gamma>0$ an $M<\infty$ such that
\begin{equation}\label{eq:B_bound}
\Pr\big(A_\rho(z,\delta,M)\big)>1-\gamma\quad\text{for every }z\in\Z^2.
\end{equation}
We further introduce the following notation for plane regions related to $\rho$:
\begin{equation*}
\begin{aligned}
H_\rho(z,\delta)&:=\big\{y\in\R^2:\rho(y-z)\le-\delta|y-z|\big\};\\
C_\rho(z,\delta)&:=\big\{y\in\R^2:|\rho(y-z)|\le\delta|y-z|\big\}.
\end{aligned}
\end{equation*}
Note that on the event $A_\rho(z,\delta,M)$ we have for all $y\in H_\rho(z,\delta)$ such that $|y-z|\ge M$ that $B_\rho(z,y)<0$. Hence, $H_\rho(z,\delta)$ corresponds to sites that are likely to be at a further distance to far out vertices along the geodesic corresponding to $\rho$ as compared to $z$.

Given $\rho_1,\rho_2,\ldots,\rho_k$ we now choose $\delta>0$ so that the cones $C_{\rho_i}(0,\delta)$, for $i=1,2,\ldots,k$, intersect only at the origin. Next, we choose $M$ large so that for all $i$
$$
\Pr\big(A_{\rho_i}(z,\delta,M)\big)>1-\eps/k.
$$
Finally, due to the choice of $\delta$ we may choose $x_1,x_2,\ldots,x_k$ so that $|x_i-x_j|\ge M$ for all $i\neq j$ and such that for each $i$ the set $H_{\rho_i}(x_i,\delta)$ contains $x_j$ for all $j\neq i$. (For instance, position the sites on a circle of large radius, in positions roughly corresponding to the directions of $\rho_1,\rho_2,\ldots,\rho_k$.) Due to these choices we will on the event $\bigcap_{i=1,2,\ldots,k}A_{\rho_i}(x_i,\delta,M)$, which occurs with probability at least $1-\eps$, have for all $i=1,2,\ldots,k$ that $B_{\rho_i}(x_i,x_j)<0$ for all $j\neq i$, as required.

It remains to show that if $|\tree_0|=\infty$ with positive probability, then it is possible to find a sequence $(x_i)_{i\ge1}$ for which $\coex(x_1,x_2,\ldots)$ occurs with probability close to one. If $|\tree_0|=\infty$ with positive probability, then it does with probability one, and $|\mathscr{C}|=\infty$ almost surely. Let $(\rho_i)_{i\ge1}$ be an increasing sequence in $\mathscr{C}$ (considered as a sequence in $[0,2\pi)$). By symmetry of $\Z^2$ we may assume that each $\rho_i$ corresponds to an angle in $(0,\pi/2)$. Fix $\eps>0$ and set $\eps_i=\eps/2^i$. We choose $\delta_1$ so that $C_{\rho_1}(0,\delta_1)$ intersect each of the lines $C_{\rho_j}(0,0)$, for $j\ge2$, only at the origin, and $M_1$ so that $\Pr(A_{\rho_1}(z,\delta_1,M_1))>1-\eps_1$. Inductively we choose $\delta_i$ so that $C_{\rho_i}(0,\delta_i)$ intersects each cone $C_{\rho_j}(0,\delta_j)$ for $j<i$ and each line $C_{\rho_j}(0,0)$ for $j>i$ only at the origin, and $M_i$ so that $\Pr(A_{\rho_i}(z,\delta_i,M_i))>1-\eps_i$. For any sequence $(x_i)_{i\ge1}$ we have
$$
\Pr\Big(\bigcap_{i\ge1}A_{\rho_i}(x_i,\delta_i,M_i)\Big)>1-\eps.
$$ 

It remains only to verify that we may choose the sequence $(x_i)_{i\ge1}$ so that for each $i\ge1$ we have $|x_i-x_j|\ge M_i$ and $x_j\in H_{\rho_i}(x_i,\delta_i)$ for all $j\neq i$. For $i\ge1$ we take $v_{i+1}\in\Z^2$ such that $|v_{i+1}|>\max\{M_1,M_2,\ldots,M_{i+1}\}$, $\rho_{i+1}(v_{i+1})>\delta_{i+1}|v_{i+1}|$ and $\rho_j(v_{i+1})<-\delta_j|v_{i+1}|$ for all $j\le i$. We note that this is possible since the sequence $(\rho_i)_{i\ge1}$ is increasing and the cone-shaped regions $C_{\rho_i}(0,\delta_i)$ and $C_{\rho_j}(0,\delta_j)$ for $i\neq j$ intersect only at the origin. Finally, take $x_1=(0,0)$, and for $i\ge1$ set $x_{i+1}=x_i+v_{i+1}$.

\section{Partial duality in higher dimensions}\label{sec:highdim}

The proof of Theorem~\ref{thm:highdim} is similar to that of Theorem~\ref{thm:duality}. So, instead of repeating all details we shall only outline the proof and indicate at what instances our current understanding of the higher dimensional case inhibits us from deriving the full duality. In the sequel we assume $d\ge2$.

The proof of the first part of the theorem is completely analogous. Suppose that
$$
\Pr\big(\coex(x_1,x_2,\ldots,x_k)\big)>0
$$
for some choice of $x_1,x_2,\ldots,x_k$ in $\Z^d$, possibly infinitely many. Then $\mathscr{N}\ge k$ almost surely, and by (Nakajima's version, which requires an exponential moment assumption, of) Corollary~\ref{cor:cardinality} we have $|\tree_0|\ge k$ almost surely.

For the second part of the argument we will need to modify slightly the approach from the two dimensional case. In the general case we do not know that every geodesic has an asymptotically linear Busemann function. However, from (the higher dimensional version of) Theorem~\ref{thm:DH} we know that if the shape has at least $k$ sides (that is, tangent hyperplanes), then, almost surely, there are $k$ geodesics in $\tree_0$ which all have asymptotically linear Busemann functions described by different linear functionals. Based on this we may repeat the proof of part~\emph{(ii)} of Theorem~\ref{thm:duality} to obtain coexistence of $k$ types with probability arbitrarily close to one.

In the case the shape has infinitely many sides, then with probability one there are infinitely many geodesics in $\tree_0$ with asymptotically linear Busemann functions, all described by different linear functionals. Let $(\rho_i)_{i\ge1}$ be a sequence of such linear functionals. Denote by $L_i$ the intersection of the hyperplane $\{x\in\R^d:\rho_i(x)=0\}$ and the $x_1x_2$-plane, i.e., the plane spanned by the first two coordinate vectors. Each $L_i$ has dimension zero, one or two, and by exploiting the symmetries of $\Z^d$ we may assume that sequence $(\rho_i)_{i\ge1}$ is chosen so that they all have dimension one. Each $L_i$ is then a line through the origin in the $x_1x_2$-plane, and by restricting to a subsequence we may assume that the sequence $(\nu_i)_{i\ge1}$ of normal vectors of these lines is monotone (considered as elements in $[0,2\pi)$). We may now proceed and select a sequence of points $(x_i)_{i\ge1}$ in the $x_1x_2$-plane in an analogous manner as in the two-dimensional case, leading to coexistence of infinitely many types with probability arbitrarily close to one.

\section*{Acknowledgements}

The author would like to express his gratitude for the years of encouragement and support received from Vladas Sidoravicius, to whom the volume in which this text will appear is dedicated to. The author is also grateful for the detailed comments received from an anonymous referee. This work was supported by the Swedish Research Council (VR) through grant 2016-04442.


\end{document}